\documentstyle{amsppt}
\magnification\magstep1
\NoRunningHeads
\pageheight{9 truein}
\pagewidth{6.3 truein}
\baselineskip=15pt
\catcode`@=11
\def\logo@{\relax}
\catcode`@=\active
\def\af #1.{\Bbb A^{#1}}
\def\au#1.{\operatorname {Aut}\,(#1)}

\def\ring#1.{\Cal O_{#1}}
\def\alb#1.{\ring .(#1)}
\def\pr #1.{\Bbb P^{#1}}
\def\prr #1.{{\Bbb P}(#1)}

\def\pic#1.{\operatorname {Pic}\,(#1)}

\def\picG#1.{{\operatorname {Pic}}^G(#1)}
\def\psl#1.{\operatorname {PSL}_{#1}}
\def\mov#1.{\operatorname {MOV}(#1)}
\def\nef#1.{\operatorname {NEF}(#1)}
\def\nd#1.{N^1(#1)}
\def\nc#1.{N_1(#1)}
\def\mc#1.{NE_1(#1)}
\def\ec#1.{\overline{NE}^1(#1)}

\def\ex#1.{\operatorname{ex}(#1)}
\def\sym#1.#2.{\operatorname{sym}_{#1}(#2)}
\def\spec#1.{\operatorname{spec}(#1)}
\def\Hom#1.#2{\operatorname{Hom}(#1,#2)}
\def\bbbq{\Bbb Q}
\def\sing#1.{\operatorname{Sing}(#1)}
\def\cox#1.{\operatorname{Cox}(#1)}

\def\bbbn{\Bbb N}
\def\bbbz{\Bbb Z}
\def\bbbc{\Bbb C}
\def\bbbcs{{\Bbb C}^{*}}

\def\sr#1.{{\Cal R}(#1)}
\def\proj#1.{\operatorname{Proj}(#1)}
\def\ac #1.{\operatorname{AMPLE}(#1)}
\def\gac #1.{C^G(#1.)}

\def\tm #1.{\operatorname{taut}(#1)}
\NoBlackBoxes

\topmatter
\title A GIT proof of W{\l}odarczyk's weighted factorization
theorem \endtitle
\author Yi Hu and Se\'an Keel \endauthor
\endtopmatter

Here we give a short proof of a 
recent result of W{\l}odarczyk, \cite{W{\l}odarczyk99a}:

\proclaim{1.0 Theorem (W{\l}odarczyk 99)} 
Let $f:X \dasharrow Y$ be a birational
rational map between smooth projective varieties in
characteristic zero. $f$ can be factored as a composition
of weighted blowups and weighted blowdowns. 
\endproclaim

For the precise definition of weighted blowup see
\cite{W{\l}odarczyk99}. Locally analytically they are 
the maps between 
quasi-smooth toric varieties given by barycentric subdivision of cones. 
See \cite{W{\l}odarzky97}.

Our goal here is a simple GIT proof of (1.0). Our main
result is:

\proclaim{1.1 Theorem} Let $f:X \rightarrow Y$ be
a birational morphism between normal $\bbbq$-factorial
varieties projective over a field.
Then there exists a normal projective variety
$Z$ with a $G= \bbbcs$ action, and two $G$-linearized
ample line bundles $L_1$, $L_2$ such that 
$Z^{ss}(L_1)//G = X$, $Z^{ss}(L_2)//G = Y$ and 
the induced rational map $X \dasharrow Y$ is $f$. Moreover
\roster
\item $Z^{ss}(L_i) = Z^{s}(L_i)$, i.e. the
GIT quotient $Z^{ss}(L_i)//G$ is geometric
\item $G$ acts freely on $Z^{ss}(L_i)$
\item The two linearizations have the same underlying
ample line bundle (i.e. they differ by a character).
\endroster 

If $X$ and $Y$ are smooth, and the characteristic is zero,
one can choose $Z$ smooth as well. 
\endproclaim

(1.0) is immediate from (1.1) and theory of wall crossings
in the variation of
GIT quotient which provides a factorization as in (1.0).
See \cite{Thaddeus96,5.6} or \cite{DolgachevHu,0.2.5}.
We need the VGIT theory (which applies to any reductive
group) only in the simplest case 
of $G = \bbbcs$. For this the
wall crossing maps are easily described using the 
Bialynicki-Birula decomposition, see e.g. \cite{Thaddeus96,1.12}.

We thank Prof. W{\l}odarczyk for sending us his preprint, which
we read before writing this note.  We first
learned of W{\l}odarczyk's work from a talk by 
Kenji Matsuki at the University of Illinois, 
where he discussed his extension, joint with Abramovich, Kalle,
and W{\l}odarczyk, of (1.0) to a factorization by blowups and 
blowdowns with smooth centers. Matsuki specifically
noted in his talk that W{\l}odarczyk's main device,
$\bbbcs$ cobordism (a pair of quotients obtained by
selecting two equivariant open subsets of a variety
with $\bbbcs$ action in a special way) admits
a GIT formulation.  The $\bbbq$-factorial 
case of (1.1) 
was known to us before.
It follows from the proof of \cite{HuKeel98,6.3}, some
of which is reproduced below -- an application of
Thaddeus's Master Space construction, \cite{Thaddeus96,3.1}. 
The smooth case follows by equivariant blowup, but we 
realized  this, and more importantly its striking
implication (1.0), only on reading W{\l}odarczyk's preprint.

We will mix 
notation of Weil divisors and line bundles. For a linearized
line
bundle $L$ and a character $v$ we let $L_v$ indicate the new
linearization obtained by twisting by $-v$. We write
$L_v^n$ for $(L_v)^{\otimes n}$.

\demo{Proof of (1.1)}
By (2.1-2.2) below, the smooth case (in characteristic
zero) follows from the general case --note if $X$ and $Y$ are
smooth then $\sing Z.$ is not semi-stable for either
linearization by (1.1.2).

Now consider the
$\bbbq$-factorial case.
Let $f: X \rightarrow Y$ be as in the statement. 
Choose an ample Cartier divisor
$D$ on $Y$. By Kodaira's lemma, see the proof of (2.2) below,
there is an effective
divisor $E$ whose support is exceptional such that
$B=f^*(D) = A + E$ with $A$ ample. 
Let $C$ be the image of
the injection ${\bbbn}^2 \rightarrow N^1(X)$ 
given by 
$(a,b) \rightarrow aA + bE$. Let $R$ be the graded ring
$$
R = \bigoplus_{v \in C} R_v = 
\bigoplus_{(a,b) \in {\Bbb N}^2} H^0(X,aA + bE).
$$
The edge generated by $B$ divides $C$ into two closed chambers:
the subcone generated by $A,B$, 
and the subcone generated by $B,E$.

\proclaim{Lemma 1} $R$ is finitely generated and after
replacing $A,B,E$ by multiples,
the canonical map
$$
R_{v_1}^{\otimes n_1} \otimes R_{v_2}^{\otimes n_2}
\rightarrow R_{v_1 + v_2} ^{\otimes n_1 + n_2} 
$$
is surjective, whenever $v_1$ and $v_2$ are in the same
chamber.
\endproclaim
\demo{Proof} Its enough to check finite generation for the 
graded subring corresponding to each chamber. By the
projection formula we reduce to the case when the two edge are
semi-ample (either $A$ and $B$ or $D$ and $\ring Y.$) 
and thus to a familiar result of Zariski. A similar argument
yields the  second statement. \qed \enddemo

Let $V = \spec R.$. 
$V$ is an instance of Thaddeus's Master Space.
Let $H = (\bbbcs)^2$. $H$ acts on $R$, with
weights $(a,b)$ on $R_{(a,b)}$. We can identify
the characters $\chi(H)$ with ${\bbbn}^2 = N^1(X)$. 
The invariants are $H^0(V,L_v)^H = R_v$ and the 
corresponding GIT quotient is 
$$
\proj {\bigoplus_{n} R_{nv}}. . \tag{1.3}
$$

\proclaim{Lemma 2} Linearizations in the interior of the same
chamber have the same semi-stable locus, which also agrees
with their stable locus. $G$ acts  freely on this locus.  The
two quotients are $X$ and $Y$. \endproclaim
\demo{Proof} 
Let $v$ be a linearization in the interior of a chamber.
A  point $h \in V$, which we identify with a $\bbbc$ algebra
map $h: R \rightarrow \bbbc$, is $v$ semi-stable
iff $h(R_{nv}) \neq 0$ for some $n > 0$.  By Lemma 1, this is true iff 
$h(R_{v'}) \neq 0$ for all $v'$ in the chamber. 
$h$ is fixed by $g =(s,t) \in H$ iff $s^a t^b =1$ for
all $v' =(a,b)$ in the chamber. Since either chamber
generates ${\bbbz}^2$ as a group, this occurs iff $g$ is
trivial. The GIT quotients are $X$ and $Y$ by (1.3) and
the projection formula. \qed \enddemo

Now if we replace $V$ by $Z = \proj R.$, with respect
to total degree, i.e. $d(a,b) = a + b$ then (1.1.1-3)
follow from Lemma 2. 
Here are the details: 

$Z$ is itself a GIT
quotient, for the diagonal subgroup
$\Delta = \bbbcs \subset H$ and line bundle
$L_v$ with $d(v) > 0$ which descends to
$\ring Z.(d(v))$. The non semi-stable locus is 
the irrelevant prime $p = \sum_{v \neq 0} R_v$,
the quotient $\pi: V \setminus p \rightarrow Z$ 
is geometric and the action on the semi-stable locus is free. 
Let $G =\bbbcs \subset H$ be a complementary 
subgroup to $\Delta$, e.g. $\bbbcs \times \{1\}$.
$H$ acts naturally on $Z$, and
$\ring Z.(1)$ has a natural $H$-linearization, such that
$\Delta$ acts trivially. Let 
$a: \chi(H) \rightarrow \chi(H/\Delta)$ be the natural
surjection. $L_v$, with its $H$ action, descends to
$\ring Z.(d(v))_{a(v)}$ ( i.e. $L_v|_{V \setminus p}$
with its $H$ action 
is the pullback of $\ring Z.(d(v))_{a(v)}$).
For any $v \neq 0$:
$p$ is $L_v$ non-stable, there are 
are natural identifications
$$
H^0(V,L_v^n)^H 
= [H^0(V,L_v^n)^{\Delta}]^{G} = 
H^0(Z,\ring Z.(d(v))_{a(v)}^n)^G,
$$
$V^{ss}(L_v) = \pi^{-1}(Z^{ss}(L_{a(v)}))$ and
$$
V^{ss}(L_v)/H = (V^{ss}(L_v)/\Delta)/G = Z^{ss}(L_{a(v)})/G. \qed 
$$
\enddemo

\remark{Remark} A slight modification of the proof of (1.1) applies
to any rational birational map. \endremark

\proclaim{2.1 Lemma} Let $f: W \rightarrow Z$ be a birational
$G$-equivariant morphism between normal projective varieties,
with $G$ reductive. Let $L$ be a $G$-linearized ample line bundle
on $W$. Suppose that $f$ is an isomorphism over
$Z^{ss}(L)$, and there exists  
an effective linearized exceptional Cartier divisor $E \subset W$ such
that $-E$ is relatively ample. 

Then there is an $n_0 > 0$ so that for
any $n>1$ and 
$M = h^*(L^{nn_0}) -E$ with
the induced linearization the following hold
for every $m>0$:
\roster
\item Every $G$-invariant section of
$M^m$ vanishes along the exceptional locus, and 
\item the canonical map
$$
H^0(W,M^m)^G @> {\otimes \sigma^{\otimes m}} >>
H^0(W,h^*(L^{mnn_0}))^G = H^0(Z,L^{mnn_0})^G
$$
is an isomorphism, where $\sigma$ is the canonical section
of $\ring.(E)$. 
\endroster
In particular $W^{ss}(M) = Z^{ss}(L)$ and
the GIT quotients are canonically identified.
\endproclaim
\demo{Proof}
We can obviously replace $L$ by a positive power.
Since 
$$
\bigoplus_{d} H^0(Z,L^{d})^G
$$
is finitely generated, after replacing $L$ by
a power 
the canonical map 
$$
\sym m. {H^0(Z,L)^G} . \rightarrow H^0(Z,L^{\otimes m})^G 
$$
is surjective for $m \geq 0$. By assumption every
invariant section of $h^*(L)$ vanishes along the 
exceptional locus, so after again replacing $L$ by a power,
invariant sections of 
$h^*(L^m)$ vanish to high order --high 
as compared with the coefficients of $mE$. The result
follows. \qed \enddemo

We expect the following is known, but we include a proof as
we do not know a reference:
\proclaim{2.2 Proposition} Let $Z$ be a normal projective 
variety in
characteristic zero, and $L$ an ample $\bbbq$-Cartier
divisor on $Z$. There exists a
resolution of singularities $f: W \rightarrow Z$, and
an effective 
Cartier divisor $E$ whose support is the full
exceptional locus, such that $-E$ is relatively ample.

If a connected uniruled group acts, this can be done equivariantly.
\endproclaim
\demo{Proof} Suppose first $Z$ is $\bbbq$-factorial, and
let $f: W \rightarrow Z$ be a projective birational map
with $W$ $\bbbq$-factorial. Let $A$ be any ample divisor
on $W$. Write
$f^*(L) -A = f^*(D) + E$ with $E$ exceptional. After adding
a multiple of $L$ to each side, we can assume $D$ is ample,
and thus absorb it into $A$. So we have
$f^*(L) = A + E$. Obviously $-E$ is relatively ample.
$E$ is effective, with support the full exceptional
locus by negativity of contraction, \cite{Koll\'areta92,2.19}.

For the general case, begin with any resolution of
singularities $f$, and arbitrary ample $A$ on $W$. Since
$f^*(L)$ is big, after replacing $L$ by a multiple,
$f^*(L) = A + D$ with $D$ effective. After
adding $f^*(L)$ to both sides, the base locus of $D$ is
contained in the exceptional locus. So we can write
$f^*(L) = A + M + E$, with $E$ effective and exceptional 
and $M$ moving,
with base locus contained in the exceptional locus. Let
$g: W \dasharrow Y$ be the rational map given by $M$. 
Let $r:W' \rightarrow W$ be a resolution of the closure of
the graph of $g$. Then  the induced $g': W' \rightarrow Y$
is regular, and $r^*(M) = M' + E'$, where  $M'$ is globally
generated, with associated map $h'$,  
and $E'$ is $r$ exceptional. 
For the composition $p:W' \rightarrow X$, we have
$p^*(L) = r^*(A) + M' + E''$, with $M'$ nef and $E''$
exceptional.
Now apply the argument of the first case to
$r^*(A)$. So we have (adjusting notation)
$p^*(L) = A + M' + E$, with $A$ ample, $M'$ nef and $E$
exceptional. $M'$ can be absorbed into $A$. Now apply
negativity of contraction as in the first case. 

Suppose connected uniruled $G$ acts, and $L$
is $G$-linearized. Note that a multiple of any line bundle on
a smooth projective $G$-variety $W$ has a linearization
by \cite{GIT,1.5}. Thus all the maps in the previous 
paragraph can be taken to be $G$-equivariant, and the
equivariant case follows. 
\qed \enddemo

\end

\end